\documentclass[10pt]{amsart}
\usepackage{amssymb}

\usepackage[colorlinks=true, pdfstartview=FitV, linkcolor=blue, citecolor=blue, urlcolor=blue]{hyperref}

\theoremstyle{plain}
\newtheorem{thm}{Theorem}[section]
\newtheorem{lemma}[thm]{Lemma}

\newtheorem{cor}[thm]{Corollary}

\numberwithin{equation}{section}
\numberwithin{figure}{section}

\theoremstyle{definition}
\newtheorem{definition}[thm]{Definition}

\newtheorem{remark}[thm]{Remark}
\newtheorem{conj}[thm]{Conjecture}

\newtheorem{qn}[thm]{Question}
\newtheorem{assumptions}[thm]{Assumptions}

\newcommand{\G}{\Gamma}
\newcommand{\Q}{\mathbb Q}

\newcommand{\Z}{\mathbb Z}
\newcommand{\R}{\mathbb R}
\newcommand{\Hyp}{\mathbb H}

\title{The large-scale geometry of right-angled Coxeter groups}
\author{Pallavi Dani}
\date{\today}
\thanks{This work was supported by a grant from the Simons Foundation (\#426932, Pallavi Dani) and  by
the Louisiana Experimental Program to Stimulate Competitive Research (EPSCoR),
funded by NSF and the Board of Regents Support Fund (NSF(2016)-LINK-108).}

\begin{document}
\begin{abstract}{This is a survey of some aspects of the large-scale geometry of right-angled Coxeter groups.  The emphasis is on recent results on their negative curvature properties, boundaries, and their quasi-isometry and commensurability classification.}
\end{abstract}
\maketitle
\section{Introduction}

Coxeter groups touch upon a number of areas of mathematics, such as representation theory, combinatorics, topology, and geometry.  These groups are generated by involutions, and have presentations of a specific form.  The connection to geometry and topology arises from the fact that the involutions in the group can be thought of as reflection symmetries of certain 
complexes.

Right-angled Coxeter groups are the simplest examples of Coxeter groups; in these, the only relations between distinct generators are commuting relations.  These have emerged as groups of vital importance in geometric group theory.  
The class of right-angled Coxeter groups exhibits a surprising variety of geometric phenomena, and at the same time, their rather elementary definition makes them easy to work with.  This combination makes them an excellent source for examples for building intuition about geometric properties, and for testing conjectures.

The book~\cite{davisbook} by Davis is an excellent source for the classical theory of Coxeter groups in  geometric group theory.  This survey focusses primarily on developments that have occurred after this book was published.    Rather than giving detailed definitions and precise statements of theorems (which can often be quite technical), in many cases we have tried to adopt a more intuitive approach.  Thus we often give the idea of a definition or the criteria in the statement of a theorem, and refer the reader to other sources for the precise details.  The aim is to give the reader a flavor of what is known about right-angled Coxeter groups.  A number of results about right-angled Coxeter groups follow from more general results about Coxeter groups.  In such cases, we have stated the right-angled version. 

The survey is organized as follows.  In Section~\ref{sec:prelim} we start with some basic definitions 
(including the definition of a right-angled Coxeter group) 
and some preliminary results that are
 used in the rest of the article.
In Section~\ref{sec:spaces}, we consider actions of 
right-angled Coxeter groups on various non-positively curved spaces.  
In particular, we give an idea of the construction of the Davis complex associated to a right-angled Coxeter group; this is a CAT(0) cube complex which admits a proper cocompact action by the group.
We also discuss hierarchical and acylindrical hyperbolicity, which are of great current interest.  In Section~\ref{sec:boundaries}, we discuss visual boundaries of CAT(0) spaces and hyperbolic spaces, and illustrate how right-angled Coxeter groups have played a role in the quest to 
explore the extent to which properties of boundaries of hyperbolic spaces hold in the CAT(0) setting.  We also discuss which topological spaces arise as boundaries of Davis complexes.  
In recent years there has been an explosion of activity related to the quasi-isometry and commensurability classification of right-angled Coxeter groups.  These results are discussed in 
Sections~\ref{sec:qi} and~\ref{sec:comm} respectively.

There is a vast body of literature on right-angled Coxeter groups, and the topics presented here are those that are closest to the expertise of the author.  Notably missing are the numerous results on 
algorithmic questions, on various notions of dimension, and on automorphisms and outer automorphisms of right-angled Coxeter groups, as well as the countless instances where right-angled Coxeter groups are used as tools or to produce examples or counterexamples for questions that arise in another area.  There is also a budding theory of random right-angled Coxeter groups,
and we mention a few results about random groups when they touch upon the topics being  discussed.

The author would like to thank Ruth Charney and Genevieve Walsh for helpful conversations,  Jason Behrstock for comments and help with Section~\ref{hhs-acyl}, and Matt Haulmark, Ivan Levcovitz, Kim Ruane, and Emily Stark for comments or corrections on a draft of this survey.   A part of this survey was written while the author was a Michler Fellow at Cornell University.  The author would like to thank the Department of Mathematics at Cornell University for its hospitality and the Association of Women in Mathematics for the Ruth I. Michler Memorial Prize.  
 
\section{Preliminaries}\label{sec:prelim}

We begin with some graph theoretic terminology and notation which is used in this survey.  
A graph is \emph{simplical} if it has no loops and there is at most one edge between each pair of vertices.  An \emph{$n$-clique} is the complete graph on $n$ vertices, denoted by $K_n$. 
Given a graph $\G$, the flag complex defined by $\G$ is the 
simplical complex which has $\G$ as its 1-skeleton, and whose $k$-simplices are in bijective correspondence with the $k$-cliques of $\G$.  If $X \subset Y$, where $X,Y$ are either graphs or complexes, then we say that $X$ separates 
$Y$ if $Y\setminus X$ has more than one component.

Given a finite simplicial graph $\G$ with vertex set $V(\G)$ and edge set $E(\G)$, the \emph{right-angled Coxeter group} $W_\G$ is the group defined by the following presentation:
$$
W_\G = \langle\, V(\G) \, \mid\, v^2=1 \text{ for all } v \in V(\G), \;[v,w]=1 \text{ if and only if } (v,w) \in E(\G) \,\rangle
$$

Throughout this survey $W_\G$ will denote a right-angled Coxeter group (even when we don't explicitly say so). The pair $(W_\G, V(G))$ is sometimes called a Coxeter system.  The \emph{nerve} of the Coxeter system is the flag complex $L=L(\G)$ defined by $\G$.  We refer to $L$ and $\G$ as the defining 
graph and defining nerve respectively.

A \emph{special subgroup} of $W_\G$ is the subgroup generated by the vertices of an induced subgraph $\Lambda$ of $\G$.  (We say that $\Lambda$ is induced if whenever two vertices of $\Lambda$ are connected by an edge of $\Gamma$, they are connected by an edge of $\Lambda$.)
The special subgroup generated by $V(\Lambda)$ is isomorphic to $W_\Lambda$.

Many algebraic and geometric features of the group $W_\G$ can be translated to properties of the defining graph.  Such properties are often called \emph{visible} or \emph{visual}.  We mention a few here pertaining to ends  and splittings of groups, that are used later in this survey. 

\begin{lemma}\label{finite} \label{finite}
$W_\G$ is 0-ended $\iff W_\G$ is finite $\iff \G$ is a clique.
\end{lemma}

Let $D_\infty$ denote the infinite dihedral group, that is $D_\infty=W_\G$ with $\G$ consisting of a pair of vertices (with no edge between them). 

\begin{thm} \label{two-ended}
 \cite[Corollary 17]{mihalik-tschantz} \cite[Theorem 8.7.3 ]{davisbook} $ W_\G$ is $2$-ended $\iff$ $W_\G$ is the product of a (possibly trivial) finite group and a $D_\infty$ 
$\iff$ $\G$ is a pair of points or the 1-skeleton of a suspension of a clique.  
\end{thm}

Stallings' Theorem says that a group has infinitely many ends if and only if it splits over a finite group.  Mihalik and Tschantz~\cite{mihalik-tschantz} showed that splittings of Coxeter groups are visible in the defining graph.  
For right-angled Coxeter groups they showed:

\begin{thm}\cite[Corollary 16]{mihalik-tschantz} \label{thm:stallings}
$W_\G$ has infinitely many ends $\iff W_G=A \ast_C B$, where $C$ is a finite special subgroup, and $A$ and  $B$ are non-trivial special subgroups $\iff \G$ has a separating 
clique.
\end{thm}

\begin{cor} \label{cor:one-ended}
$W_\G$ is one-ended $\iff$ $\G$ is not equal to a clique, is connected and has no separating cliques.  
\end{cor}

Moreover, any right-angled Coxeter group has a particularly nice graph of groups decomposition:

\begin{thm}\cite[Theorem 8.8.2 ]{davisbook}\cite[Corollary 18]{mihalik-tschantz}\label{one-ended}
A right-angled Coxeter group with infinitely many ends splits as a tree of groups where each vertex group is a 0- or 1- ended special subgroup, and each edge group is a finite special subgroup. 
\end{thm}

We will also be interested in splittings over 2-ended subgroups.  It follows from the Main Theorem 
of~\cite{mihalik-tschantz} that: 

\begin{thm}\cite{mihalik-tschantz}\label{two-ended-splittings}
$W_\G$ splits over a 2-ended group $\iff$ $\G$ has a separating subgraph which is either a pair of points or the 1-skeleton of a suspension of a clique.  
\end{thm}

Finally, in a different direction, we mention that the virtual cohomological dimension (vcd) of a right-angled Coxeter group can also be computed directly from its nerve, by the following formula due to Davis.  
\begin{thm}\cite[Corollary 8.5.5]{davisbook}\label{vcd} 
Given a right-angled Coxeter group $W_L$ with nerve $L$,
$$
\mathrm{vcd} (W_L) = \mathrm{max} \{ n \, |\, \overline{H}^{n-1}(L-\sigma) \neq 0, \text{ for some simplex } \sigma \text{ of } L\}.
$$
\end{thm}

\medskip
\section{Actions of right-angled Coxeter groups on non-positively curved spaces}
\label{sec:spaces}

Exploiting properties of negatively or non-positively curved spaces to draw conclusions about groups 
which act on them geometrically is a central theme in geometric group theory.  
Recall that a group action is \emph{geometric} if it is a proper cocompact action by isometries. 
In this section we explore geometric actions of right-angled Coxeter groups on non-positively curved spaces.  
We assume that the right-angled Coxeter group is not finite.  By Lemma~\ref{finite} this is equivalent to the condition that the defining graph has at least one pair of vertices which are not connected by an edge. 

\subsection{CAT(0) cube complexes, virtual specialness}

All Coxeter groups (not just the right-angled ones) are CAT(0) groups: 
every Coxeter group acts properly and cocompactly on its~\emph{Davis complex}, which is CAT(0).  
When the Coxeter group is right-angled, the Davis complex can be metrized to have the structure of a CAT(0) cube complex.  
The precise definition of the Davis complex and the proof that it is CAT(0) can be found 
in~\cite{davisbook},  in Chapters 7 and 12 respectively.  Here we give an idea of the 
construction in the right-angled case as described in Chapter~1 of~\cite{davisbook}.  It is reminiscent of the construction of the Salvetti complex associated to a right-angled Artin group.  

\begin{definition}[Davis complex $\Sigma_\G$]\label{davis-cplx}
Given a graph $\G$ with $n$ vertices, let $L$ be the flag complex on $\G$.  
We start by constructing a certain subcomplex $P_L$ of $C=[-1,1]^n\subset \R^n$.  If we fix an order on 
the vertex set $V(\G)$ of $\G$, so that the $i$th generator of $W_\G$ corresponds to the 
$i$th factor of $C$, then each face of $C$ determines a unique subset of $V(\G)$.  Define $P_L$ 
to be the subcomplex of $C$ which contains the face $F$ of $C$ if and only if the subset of $V(\G)$ corresponding to $F$ is a clique.  Then the Davis complex $\Sigma_\G$ of $W_\G$ is 
the universal cover of $P_L$.  

The group $(\Z_2)^n$ acts on $C$ (and hence $P_L$) 
by reflections in hyperplanes of $\R^n$ passing through the origin.  This 
 induces a finer cube complex structure on $P_L$, which lifts to a cube complex structure on 
 $\Sigma_\G$, and the $\Z_2^n$-action on $P_L$ lifts to an action of $W_\G$ on $\Sigma_\G$ by cubical isometries.  The quotient of $\Sigma_\G$ by the $W_\G$ action is equal to the quotient of 
 $P_L$ by the $\Z_2^n$-action.  One can verify that $\Sigma_\G$ is CAT(0) by checking the Gromov link condition.   
 \end{definition}
 
 Thus right-angled Coxeter groups are CAT(0) cube complex groups.  
An important class of CAT(0) cube complex groups consists of the \emph{virtually special} ones.  
These are groups which have a finite-index subgroup that is the fundamental group of a 
\emph{special cube complex}.  Haglund and Wise introduced special cube complexes in~\cite{haglund-wise1} in terms of the absence of certain hyperplane pathologies.  One can check directly that these pathologies are absent in $P_L$, and as a result, right-angled Coxeter groups are virtually (compact) special.  

We remark that even the non-right-angled  Coxeter groups act on CAT(0) cube complexes.  
For any Coxeter group, Niblo and Reeves~\cite{niblo-reeves} defined a CAT(0)  cube complex 
on which the Coxeter group acts properly by cubical isometries, but in general the action is not cocompact.  For right-angled Coxeter groups $W_\G$,
these complexes agree with the Davis complex $\Sigma_\G$ defined above.  
In~\cite{haglund-wise2}, Haglund and Wise used these complexes to show that all Coxeter groups are virtually special.  A long standing open question about general Coxeter groups is whether they are biautomatic.  This would be implied by showing that they are virtually compact special.

\subsection{Hyperbolicity}
We now explore which right-angled Coxeter groups exhibit various degrees of negative curvature, starting with the strongest condition: that of being (Gromov) hyperbolic. 
Moussong characterized the hyperbolic right-angled Coxeter groups in terms of their defining graphs:
\begin{thm}\cite{moussong}\label{moussong}
A right-angled Coxeter group $W_\G$ is hyperbolic if and only if $\G$ has no induced squares.
\end{thm}

A very special class of hyperbolic groups consists of those that act geometrically on the hyperbolic space $\mathbb H^n$ of some dimension $n$.   A detailed account of the classical theory of Coxeter groups acting on the three constant curvature geometries ($S^n$, $\mathbb{E}^n$ and $\mathbb H^n)$ can be found in~\cite[Chapter 6]{davisbook}. 
We remark here that by work of Vinberg, 
if a right-angled Coxeter group acts geometrically on $\Hyp^n$, then $n\le 4$ 
(see~\cite[Corollary~6.11.7]{davisbook}).  More generally, Vinberg~\cite{vinberg} showed that $n \le 29$ for an 
arbitrary  (non-right-angled) Coxeter group, but examples are known only up to $n=8$.  
Januszkiewitz and \'{S}wi\c{a}tkowski~\cite{janus-swiat} show that if one passes to the 
class of hyperbolic Coxeter groups which are Poincare\'e duality groups of dimension~$n$, then one still has a bound on the dimension (of $n=4$ in the right-angled case, and $n=61$ in general).

Vinberg's bound on the dimension of an $\Hyp^n$ admitting a geometric Coxeter group action led Moussong to conjecture in~\cite{moussong} that there is also a universal bound on the virtual cohomological dimension of a 
Gromov hyperbolic right-angled Coxeter group, and this question was also raised in~\cite{gromov}.  
Januszkiewicz and \'{S}wi\c{a}tkowski disproved this 
in~\cite{janus-swiat} by constructing examples of hyperbolic right-angled Coxeter groups with vcd$=n$, for any $n$.   To do this they used the theory of complexes of groups to construct flag simplicial complexes $L_n$ for any $n$, which satisfy the conditions in 
Davis' formula for the vcd$=n$ (Theorem~\ref{vcd} above) and Moussong's characterization of hyperbolicity
(Theorem~\ref{moussong} above). 

Subsequently there have been variations of these constructions in~\cite{haglund, janus-swiat2, abj+}, which all use complexes of groups machinery.   In~\cite{osajda} Osajda gave a new construction of hyperbolic 
 right-angled Coxeter groups of arbitrarily high vcd.  His defining complexes are 
 constructed by an elementary inductive procedure that forces the dimension of the complex to go up each time, and does not use complexes of groups.   Osajda points out that the previously 
 known constructions all yielded examples which are systolic, a condition which forces these groups to be ``asymptotically two-dimensional'' in a certain sense.  He used the procedure in~\cite{osajda} 
to produce examples which are not systolic.  

\subsubsection{Random right-angled Coxeter groups and hyperbolicity}\label{random}
Gromov~\cite{gromov} introduced the theory of random groups in terms of random presentations, and showed that the random group (in the few relators model) is hyperbolic.  In the density model, 
in which the number of relations is allowed to grow with the length of the relations, one has
a ``phase transition'': when the density is $<1/2$, the random group is infinite and hyperbolic, while when the density is $>1/2$, the random group is finite.  

The more natural way to think about random right-angled Coxeter groups is as those defined by random graphs, for which there is a well-developed theory.  In the Erd\H{o}s--R\'enyi model,
the probability space $G(n,p(n))$ consists of graphs on $n$ vertices, 
with the probability measure coming from 
independently choosing edges between pairs of vertices with probability $p(n)$.  A property holds asymptotically almost surely (a.a.s.) if the probabililty that it holds for graphs in $G(n, p(n))$ tends to 1 as $n$ tends to infinity.  Charney and Farber investigated the hyperbolicity of random right-angled Coxeter groups with this model, and showed: 
\begin{thm}\cite{charney-farber}
The right-angled Coxeter group corresponding to  a graph in $G(n, p(n))$ is hyperbolic a.a.s.~ if either 
$p(n)(1 - n)^2 \to 0$ or $p(n)n \to 0$.  On the other hand, if $p(n)n \to \infty$ and $(1-p(n))n^2 \to \infty$ then the random right-angled Coxeter group is a.a.s.~not hyperbolic. 
\end{thm}

\subsection{Relative hyperbolicity}
Let $G$ be a finitely presented group and let $\mathcal P$ be a finite collection of proper subgroups of $G$.  
Recall that $G$ is said to be \emph{hyperbolic relative to $\mathcal P$} if the \emph{coned-off Cayley graph} (obtained by collapsing the cosets of $\mathcal P$ to finite diameter sets in a Cayley graph of $G$) is $\delta$-hyperbolic, and if the collection~$\mathcal P$ satisfies the \emph{bounded coset penetration property}, which controls the behavior of quasi-geodesics in the coned-off Cayley graph.  
A number of definitions of relative hyperbolicity can be found in the literature.  
See~\cite{hruska} for the details of these definitions, and a proof of their equivalence. 
The relatively hyperbolic structure 
$(G, \mathcal P)$ is said to be \emph{minimal}, if  given any other relatively hyperbolic structure $(G, \mathcal Q) $
on~$G$, every group in $\mathcal P$ is conjugate to some group in $\mathcal Q$.  

Caprace obtained a criterion for the relative hyperbolicity of general Coxeter groups in~\cite[Theorem A$'$]{caprace-erratum}.  (Note that this corrects~\cite[Theorem A]{caprace}.)  The following is the statement of Caprace's characterization applied to the special case of right-angled Coxeter groups, and stated in terms of the defining graph, rather than in terms of Coxeter systems.  To translate from the statement in~\cite[Theorem A$'$]{caprace-erratum} we use the fact that the only irreducible affine 
right-angled Coxeter group is $W_\Gamma=D_\infty$, for which $\G$ consists of a pair of vertices not connected by an edge, which makes the first condition simpler than the general case. 
In the statement below, define a graph to be irreducible if it does not split as a join, and for any graph 
$J \subset \G$, define $J^{\perp}$ to be the subgraph of $\G$ induced by the generators in $V(\G) \setminus V(J)$ which commute with every element of $V(J)$ (or equivalently, are connected by an edge of $\G$ to every element of $V(J)$).

\begin{thm}\cite[Theorem A$'$]{caprace-erratum}\label{caprace}
Let $W_\G$ be a right-angled Coxeter group, and let $\mathcal J$ be a collection of proper 
induced subgraphs of $\G$. Then $W$ is hyperbolic relative to $\{W_J \mid J \in \mathcal J\}$ if and only if  the following conditions 
hold. 
\begin{enumerate}
\item For any pair of irreducible non-clique subgraphs $L_1$ and $L_2$ of $\G$ such that $\G$ contains the join of $V(L_1)$ and 
$V(L_2)$ (i.e.~each vertex of $L_1$ is connected by an edge of $\G$ to each vertex of $L_2$), there exists 
$J \in \mathcal J$ such that $L_1, L_2 \subset J$. 
\item For all $J_1, J_2 \in \mathcal J$ with $J_1\neq J_2$, the intersection $J_1 \cap  J_2$ is a (possibly empty) clique.
\item For each $J \in \mathcal J$ and each irreducible non-clique subgraph $L \subset J$, we have 
$L^{\perp} \subseteq J$.
\end{enumerate}
\end{thm}

In~\cite[Corollary D]{caprace}, Caprace additionally characterizes exactly those Coxeter groups which are relatively hyperbolic with isolated flats.  

The notion of \emph{thickness} was introduced 
in~\cite{behr-dru-mos} 
as a geometric obstruction to relative hyperbolicity.  
An idea of the definitions of thickness and strong algebraic thickness  
is given in Section~\ref{divthick}.  
Behrstock, Caprace, Hagen,  and Sisto showed in~\cite[Apendix A]{behr-hag-sis}
that for Coxeter groups, there is in fact a dichotomy: they are either thick, or relatively hyperbolic relative to thick groups.  
In the right-angled case, their theorem says:

\begin{thm}\cite{behr-hag-sis}\label{relhypthick}
Let $\mathcal T$ be the class consisting of the finite simplicial graphs $\Lambda$ such that 
$W_\Lambda$ is strongly algebraically thick. Then for any finite simplicial graph $\Gamma$, either: $\G \in  \mathcal T$, or there exists a collection $\mathcal J$ of induced subgraphs of $\G$ such that $\mathcal J \subset \mathcal T$ and $W_\G$ is hyperbolic relative to the collection $\{W_J | J \in \mathcal J\}$.  In the latter case, this relatively hyperbolic structure is minimal.\end{thm}

Moreover, Behrstock--Hagen--Sisto provide a characterization of the graphs in $\mathcal T$ in terms of an inductively defined 
property of the defining graph~\cite[Theorem II]{behr-hag-sis}.  Using this characterization, they obtain a polynomial time algorithm to decide if a right-angled Coxeter group given by a graph is relatively hyperbolic or thick~\cite[Theorem III]{behr-hag-sis}.

Thickness and relative hyperbolicity are both quasi-isometry invariants, and as a result, this 
Theorem~\ref{relhypthick} above provides 
a big step in the quasi-isometry classification of right-angled Coxeter groups.  This is discussed further in 
Sections~\ref{divthick} and~\ref{relhyp}.   In~\cite{behr-hag-sis,bfhs}, the authors study when random right-angled Coxeter groups are relatively hyperbolic or thick.  These results are discussed in Section~\ref{random-rh}.

\subsection{Hierarchical and acylindrical hyperbolicity}\label{hhs-acyl}
Behrstock, Hagen and Sisto introduced 
\emph{hierarchically hyperbolic spaces} in~\cite{bhs-hhs} 
as a common framework for studying mapping class groups and cubical groups. 
A recent survey by Sisto~\cite{sisto-hhs} provides an excellent introduction to the area.  
As the name suggests, a hierarchically hyperbolic space can be thought of as being built up inductively, starting from a base level consisting of hyperbolic spaces.  At the top level, the hierarchically hyperbolic space is weakly hyperbolic relative to a collection of regions called standard product regions, meaning that the coned-off space is hyperbolic but the bounded coset penetration property may not be satisfied.  The standard product regions are themselves hierarchically hyperbolic spaces of lower ``complexity''.  This setup is suited to adapting the powerful machinery developed by Masur and Minsky to study mapping class groups, and allows for inductive arguments.  

Behrstock--Hagen--Sisto introduced 
the notion of a  \emph{factor system} for a cube complex, which is a collection of convex subcomplexes that is uniformly locally finite, contains all combinatorial hyperplanes, and satisfies a certain projection criterion.  They showed in~\cite[Theorem G]{bhs-hhs} that any cube complex which has a factor system is a hierarchically hyperbolic space, and 
in~\cite[Proposition~B]{bhs-hhs} that the universal cover of any compact special cube complex has a factor system.  
It follows that right-angled Coxeter groups are hierarchically hyperbolic.  

Another class of groups of great current interest is that of \emph{acylindrically hyperbolic groups}.  In~\cite{osin}, Osin showed that the class of groups which admit a non-elementary, acylindrical action on a hyperbolic space coincides with a number of different generalizations of relative hyperbolicity in the literature, and initiated a study of these groups.  It is shown 
in~\cite{bhs-hhs} that under some minimal conditions (which hold for all right-angled Coxeter groups except those which are direct products), 
hierarchically hyperbolic groups admit non-elementary acylindrical actions.

In~\cite{abb-behr-dur} Abbott, Behrstock and Durham show that any hierarchically hyperbolic group (in particular, a right-angled Coxeter group which is not a product) admits a \emph{universal acylindrical action}, i.e., an acylindrical action such that if a group element acts loxodromically in some acylindrical action on some hyperbolic space, then it acts loxodromically in this action.  Moreover, they show that there exists a \emph{largest acylindrical action}, i.e., an action which is comparable to, and larger than all other 
acylindrical actions, in a particular partial order on cobounded acylindrical actions.  In particular, a largest action is universal and unique.

The areas of hierarchically and acylindrically hyperbolic groups are evolving quickly, and the resulting theorems should have interesting implications for right-angled Coxeter groups.

\section{Boundaries of right-angled Coxeter groups}\label{sec:boundaries}

Let $X$ be a metric space.  The set of boundary points $\partial X$ of $X$ consists of all equivalence classes of geodesic rays in $X$, where two rays are equivalent if and only if they are asymptotic. When $X$ is a CAT(0) space, there is a natural topology on 
$X\sqcup \partial X$, called the cone topology.  Roughly speaking, two geodesic rays  are close in the cone topology if they fellow travel for a long distance.  The \emph{visual boundary} of $X$ is the set $\partial X$ endowed with the cone topology.  When $X$ is hyperbolic in addition, this agrees with the Gromov boundary.  See~\cite[Chapter II.8]{bridson-haef} for detailed definitions.

As described in Section~\ref{hyperbolic} below, the boundaries of hyperbolic spaces which admit a geometric action by a group behave extremely well.  They have some strong topological restrictions, and some of these can be used to draw conclusions about the algebraic structure of the group.  
There has been a lot of work in the past two decades to understand the extent to which these properties carry over to the non-positive curvature setting.  While some of these questions are hard if one considers all CAT(0) spaces, right-angled Coxeter groups can provide some insight into what might be expected in general. 

\bigskip

\subsection{Boundaries of hyperbolic groups}~\label{hyperbolic}
The homeomorphism type of the visual boundary of a hyperbolic space is a quasi-isometry invariant.  More precisely:
\begin{thm}\label{hyp-quasi}
Let $X$ and $X'$ be hyperbolic spaces.  If $f: X \to X'$ is a quasi-isometry then $f$ extends to a homeomorphism 
from $\partial X$ to $\partial X'$.  

Moreover, if $G$ is a group acting geometrically on two hyperbolic spaces $X$ and $X'$, then the $G$-action extends to the boundaries $\partial X$ and $\partial X'$, and
the natural quasi-isometry between $X$ and $X'$ (through $G$) induces a $G$-equivariant homeomorphism 
$\partial X \to \partial X'$.  
\end{thm}
In particular, there is a well-defined notion of the boundary $\partial G$ of a hyperbolic group $G$ (up to homeomorphism).  

Hyperbolicity imposes some additional topological restrictions on the boundary.  
By work of Bestvina--Mess~\cite{bestvina-mess} and Swarup~\cite{swarup} one has that 
if $G$ is a one-ended hyperbolic group, then $\partial G$ has no global cut points and is locally 
connected.  Topological features of the boundary are sometimes manifested in algebraic features of the group.  For instance, Bowditch~\cite{bowditch} shows that 
if $G$ is a one-ended hyperbolic group which is not cocompact Fuchsian, then local cut points in $\partial G$ correspond to splittings of $G$ over two-ended subgroups.   

Kapovich--Kleiner obtain the following characterization of 1-dimensional boundaries of hyperbolic groups.

\begin{thm}\cite{kapovich-kleiner}\label{thm:kk}
If $G$ is a hyperbolic group with 1-dimensional boundary which does not split over a finite or 2-ended group,
then $\partial G$ is either a circle,  
a Sierpinski carpet or a Menger curve.\end{thm}
\noindent
The proof of Theorem~\ref{thm:kk} uses the 
topological properties of hyperbolic boundaries in the previous paragraph together with the following topological characterizations of the Sierpinski carpet and the Menger curve.

\begin{thm}\label{sierp-meng}\cite{whyburn, anderson1, anderson2}
Any compact, 1-dimensional, planar, connected,
locally connected space with no local cut points is a Sierpinski carpet. 

A compact, metrizable,
connected, locally connected, 1-dimensional space is a Menger curve provided it has
no local cut points, and no nonempty open subset is planar.  
\end{thm}

The survey~\cite{benakli-kapovich} provides a far more detailed account of boundaries of hyperbolic groups.

\subsection{Boundaries in the CAT(0) setting}\label{cat0boundary}

In contrast with the hyperbolic case, the boundaries of CAT(0) spaces which admit geometric 
group actions are quite badly behaved.  By now  it is well known that even the weakest 
generalization of Theorem~\ref{hyp-quasi} fails:
Croke--Kleiner~\cite{croke-kleiner} give an example of a group acting geometrically on two 
CAT(0) spaces whose boundaries are not homeomorphic.   In fact, 
Wilson~\cite{wilson} showed that the group from~\cite{croke-kleiner} admits uncountably many distinct boundaries. 

If one considers geometric actions of a group $G$ on two CAT(0) spaces $X$ and $Y$ which are known to have homeomorphic boundaries, one might ask whether these boundaries are equivariantly homeomorphic, and moreover, whether the natural quasi-isometry between $X$ and $Y$ (through $G$) extends to a homeomorphism of their boundaries, as in Theorem~\ref{hyp-quasi}.   Both of these things turn out to be false in general.  Croke--Kleiner~\cite{croke-kleiner2} and Buyalo~\cite{buyalo}  gave examples of such $X, Y$ and $G$ with the property that $\partial X$ and $\partial Y$ are homeomorphic but not $G$-equivariantly homeomorphic.  
Bowers--Ruane \cite{bowers-ruane} considered a CAT(0) space $X$ of the form $\mathrm{Tree}\times \mathbb{R}$ with two different actions
by $G(=F_2 \times \Z)$.  The first was a simple product action and the other was a  ``twisted action'' obtained by precomposing with an element of the automorphism group of $G$.  
They showed that  there exists a $G$-equivariant homeomorphism $\partial X \to \partial X$.  However, 
the 
natural quasi-isometry from $X$ to itself induced by these actions does not extend to a continuous map on the boundary.

The groups constructed in~\cite{bowers-ruane, buyalo, croke-kleiner, croke-kleiner2} are not right-angled Coxeter groups, but such phenomena exist in the right-angled Coxeter group setting as well.  
By considering a right-angled Coxeter group commensurable to $F_2 \times \Z$, and using the same method as 
in~\cite{bowers-ruane}, Yamagata~\cite{yamagata} constructed a right-angled Coxeter version of the Bowers-Ruane example.  Thus one has:

\begin{thm}\cite{yamagata}
There exists a right-angled Coxeter group which admits two different geometric actions on a space $X$, such that the natural quasi-isometry of $X$ induced by the two group actions does not extend to a continuous map $\partial X  \to \partial X$. 
\end{thm}

In~\cite{qing1} and~\cite{qing2} Qing studied the actions of right-angled Coxeter groups on the \emph{Croke-Kleiner spaces} from~\cite{croke-kleiner}.  The group in the Croke-Kleiner example is the  right-angled 
Artin group whose defining graph is the path of length four.  The corresponding Salvetti complex consists of three tori in a chain, so that the middle torus is glued to each of the two other tori along a simple closed curve, and the two simple closed curves in the middle torus intersect exactly once.  (See~\cite[Fig 3]{qing1} for a picture.) 
One then obtains an uncountable family of CAT(0) spaces by allowing the angle between the two simple closed curves to be anything in $(0, \pi/2]$ and passing to the universal cover.  Each such space still admits a
geometric action by the same right-angled Artin group, but Croke and Kleiner showed in~\cite{croke-kleiner} that the resulting boundary when the angle is $\pi/2$ is not homeomorphic to the resulting boundary when the angle is not equal to $\pi/2$.  Wilson showed in~\cite{wilson} that in fact, no two angles result in spaces with homeomorphic boundaries.  

In~\cite{qing2}, Qing shows that right-angled Coxeter groups are more ``geometrically rigid'' than right-angled Artin groups, in the following sense:
if a right-angled Coxeter group acts geometrically on one of the Croke-Kleiner spaces above, then the angle must equal $\pi/2$.  On the other hand, in~\cite{qing1} she considers the  universal covers 
of the spaces obtained by keeping the angle between the two curves 
equal to $\pi/2$ but varying the lengths of the simple closed curves in the central torus (so that the tori are built out of rectangles instead of squares).  
These spaces are still CAT(0). 
She shows:
\begin{thm}\cite{qing1}
There exists a right-angled Coxeter group which acts geometrically on all the Croke-Kleiner spaces obtained by varying the lengths of the two simple closed curves in the central torus, but the boundaries of these spaces are not all equivariantly homeomorphic.  
\end{thm}

It is not known whether the boundaries of the spaces from~\cite{qing1} are homeomorphic, so the following question is still open: 

\begin{qn}
Is there a right-angled Coxeter group with acts geometrically on two spaces whose boundaries are not homeomorphic?
\end{qn}

\subsection{Local connectedness}
The results in Section~\ref{cat0boundary} show that the boundary of a CAT(0) group is not well-defined and is not a quasi-isometry invariant in general. 
Nevertheless, 
one could still ask whether there are topological features which are shared by all the boundaries associated with a particular group.  
One such property is local connectedness.  

As mentioned in Section~\ref{hyperbolic}, the boundary of a one-ended hyperbolic group is always locally connected.  
The authors of~\cite{camp-mihalik} point out that 
although the Croke--Kleiner group described in Section~\ref{cat0boundary} admits uncountably many non-homeomorphic boundaries, they are all non-locally connected. 
In fact, the following is still open:
\begin{qn}
Is there a CAT(0) group which acts on two different spaces, such that one has locally connected boundary and the other has non-locally connected boundary?
\end{qn}

There are easy examples of 1-ended groups acting geometrically on CAT(0) spaces whose boundaries are not locally connected.  For example, 
consider the group $F_2 \times Z$ acting on $X= T\times \R$, where $T$ is the Cayley graph of $F_2$.  
Then $\partial X$ is the suspension of a Cantor set, and therefore it fails to be locally connected 
everywhere except at the two suspension points. 
Observe that $F_2 \times \Z$  can also be written as an amalgamated product of the form $\Z^2 \ast_\Z \ast \Z^2$.  
Mihalik--Ruane~\cite{mihalik-ruane} generalized this example to show that whenever a group admits a particular type of ``geometric'' splitting as an amalgamated product, the boundary of any space on which the group acts geometrically is non-locally connected.  (See~\cite[Main Theorem]{mihalik-ruane} or~\cite[Theorem~3.3]{mrt} for a more general version.) They also applied this to the case of right-angled Coxeter groups and identified a sufficient condition on the defining graph which guarantees non-locally connected boundary.

In~\cite{mrt}, the authors asked whether there is a converse to the Mihalik-Ruane geometric splitting theorem.  
\begin{qn}\cite{mrt}
Given a CAT(0) space $X$, 
can the non-local connectedness of the boundary be completely characterized in terms of the existence of a ``geometric'' splitting of groups which act on $X$ properly and cocompactly?
\end{qn}
In this generality this question remains open but progress has been made in some cases, as described below.

\subsubsection{Isolated flats}\label{isolated flats}

In a recent paper, Hruska--Ruane~\cite{hruska-ruane} show that given a group $G$ acting geometrically on a one-ended CAT(0) space $X$ with isolated flats, if $\partial X$ is not locally connected, then $G$ admits a Mihalik--Ruane splitting.  Thus they obtain 
a complete characterization of non-locally connected boundaries in the isolated flats setting.  

\subsubsection{Right-angled Coxeter groups}
In~\cite{mrt}, Mihalik, Ruane and Tschantz sought a converse to the Mihalik--Ruane theorem in the case of one-ended right-angled Coxeter groups.  They defined two combinatorial objects associated to the defining graph of a right-angled Coxeter group: \emph{product separators} and \emph{virtual factor separators}.  They point out that if $W$ splits as a visual direct product, then $W$ has non-locally connected boundary if and only if one of the factor groups does.  
In all other cases, they showed the following.

\medskip
\begin{thm}~\cite{mrt}
Let  $W_\G$ be a one-ended subgroup which does not split visually as a direct product.
If $\G$ contains a virtual factor separator, then $W_\G$ has non-locally connected boundary.  On the other hand, if $\G$ contains no product separator and no virtual factor separator, then $W_\G$ has locally connected boundary.  
\end{thm}
Camp and Mihalik~\cite{camp-mihalik} promoted this to a full converse to the Mihalik-Ruane theorem under the additional hypothesis that $W$ contains no
$\Z^3$ subgroups.  

\begin{thm}\cite{camp-mihalik} 
If $W_\G$ is a one-ended right-angled Coxeter group which has no 
visual subgroup isomorphic to $(\Z_2 \ast \Z_2)^3$ and does not split visually as a direct product, then $W_\G$ has locally connected boundary if and only if $\G$ has no virtual factor separator.  
\end{thm}

In~\cite[Corollary D]{caprace}, Caprace characterizes exactly when a right-angled Coxeter group is CAT(0) with isolated flats.  Combining this with Hruska and Ruane's result in Section~\ref{isolated flats}
will provide further examples of right-angled Coxeter groups for which the converse  to the Mihalik--Ruane theorem holds.  
However, since there exist non-hyperbolic right-angled Coxeter groups with 2-dimensional flats that are not isolated, one still does not have the full converse 
for right-angled Coxeter groups.

\subsection{Topological spaces which arise as boundaries}

Kapovich and Kleiner's result (Theorem~\ref{thm:kk} above) says that the only spaces which arise in dimension 1 as boundaries of ``indecomposable'' hyperbolic groups are $S^1$, the Sierpinski carpet, and the Menger curve.  This  prompts some natural questions.  
Can one characterize which hyperbolic groups have one of these spaces as their boundary?  Which spaces arise as boundaries of hyperbolic groups in dimension higher than 1?  More generally, which spaces arise as boundaries of CAT(0) spaces which admit a geometric group action?

The results of \cite{tukia, gabai, casson-jungreis} show that if $G$ is a hyperbolic group with $\partial G =S^1$ such that $G$ acts as a convergence group on $\partial G$, then $G$ is cocompact Fuchsian.  Champetier~\cite{champetier} showed that the generic finitely presented group has 
Menger curve boundary.  In light of this Kapovich--Kleiner remark in~\cite{kapovich-kleiner} that 
it is ``probably impossible to classify hyperbolic groups whose boundaries are homeomorphic
to the Menger curve''.  On the other hand, characterizing groups with Sierpinksi carpet boundary is more approachable, though still challenging (and still open).   Kapovich--Kleiner conjectured that if 
$G$ is a hyperbolic group with Sierpinski carpet boundary, then G acts
geometrically on a convex subset of $\Hyp^3$ with non-empty totally geodesic boundary.

For right-angled Coxeter groups one can say more, including in the non-hyperbolic case.  
Since the boundary is not necessarily well-defined in the non-hyperbolic case, we use the following convention in this section:
\begin{remark}
In the discussion below, when the right-angled Coxter group is not hyperbolic, ``the boundary'' refers to the boundary of the Davis complex.  
\end{remark}
Define a simplicial complex to be \emph{inseparable} if it is connected, has no separating subcomplex 
which is a simplex, a pair of points, or a suspension of a simplex.
\'{S}wi\c{a}tkowski used this property to characterize when a right-angled Coxeter group with planar nerve has 
Sierpinski carpet boundary:  
 \begin{thm}\cite[Corollary 1.4]{swiat2}\label{racg-sierp}
If $W$ is a right-angled Coxeter group with planar nerve $L$, then 
$\partial W$ is a Sierpinski carpet if and only if $L$ is inseparable and distinct from a simplex and from a triangulation of $S^2$.  
 \end{thm}
This generalizes~\cite[Theorem 1]{swiat1}, which characterizes the hyperbolic right-angled Coxeter groups with planar nerves whose boundary is homeomorphic to the Sierpinski carpet.  
Note that by Theorems~\ref{thm:stallings} and~\ref{two-ended-splittings} above, the inseparability of the nerve is equivalent to the fact that the group does not split over a finite or 2-ended group.  In the hyperbolic setting (in~\cite{swiat1}), \'{S}wi\c{a}tkowski uses this fact, together with the topological properties of hyperbolic boundaries, 
the Kapovich--Kleiner theorem (Theorem~\ref{thm:kk} above), and the topological characterization of the 
Sierpinski carpet (Theorem~\ref{sierp-meng} above), 
 to reduce the proof of Theorem~\ref{racg-sierp} to showing that when the nerve is planar and not a simplex, the boundary of the group is planar and 1-dimensional. For the latter property, Theorem~\ref{vcd} above is used to conclude that the vcd of the group is~2.  

In~\cite{swiat1}, \'{S}wi\c{a}tkowski also provides a conjecture for when an arbitrary (i.e.~not necessarily hyperbolic and not-necessarily right-angled) Coxeter group with planar nerve has Sierpinski carpet boundary.  
Haulmark~\cite[Theorem 1.5]{haulmark} proves \'{S}wi\c{a}tkowski's conjecture for Coxeter groups with isolated flats. 
This is an application of the Main Theorem of Haulmark's paper~\cite{haulmark}, which is a version of the Kapovich--Kleiner characterization of 1-dimensional boundaries for groups acting on CAT(0) spaces with isolated flats.   Note that  \'{S}wi\c{a}tkowski's result (Theorem~\ref{racg-sierp} above) addresses the right-angled case of the conjecture.  Haulmark's result includes non-right-angled Coxeter groups.

Planarity of the nerve is an assumption in \'{S}wi\c{a}tkowski's results characterizing the right-angled Coxeter groups with Sierpinski carpet boundary.  He speculates in~\cite{swiat1}
that up to taking a product with a finite group, these may be exactly the hyperbolic groups which have Sierpinski carpet boundary.  
In upcoming work~\cite{dhw}, the author, Haulmark, and 
Walsh prove a complementary result, showing that if 
$\G$ is triangle-free and non-planar, and the flag complex $L(\G)$ is inseparable, then the boundary of $W_\G$ is non-planar. 
Combining this with the results of Kapovich--Kleiner~\cite{kapovich-kleiner} and 
Haulmark~\cite{haulmark}, and the characterization of the Menger curve in Theorem~\ref{sierp-meng} above, it follows that if 
$\G$ is triangle-free and non-planar, the flag complex $L(\G)$ is inseparable, and
$W_\G$ is either hyperbolic or has isolated flats, then the boundary of $W_\G$ is a Menger curve.  
Haulmark-Nguyen-Tran~\cite{haul-ngu-tra} also give examples of non-hyperbolic right-angled Coxeter groups with Menger curve boundary.

In dimensions $n>1$, there are still very few explicit topological spaces of dimension $n$ which are 
known to be boundaries of hyperbolic groups.  Some known examples are spheres and Sierpinksi compacta (which are boundaries of fundamental groups of hyperbolic manifolds), and the Menger compacta $M_{2,5}$ and $M_{3,7}$~\cite{dymara-osajda} (which are discrete cocompact automorphism groups of right-angled hyperbolic buildings).  In addition, Dranishnikov~\cite{dran} constructed hyperbolic right-angled Coxeter groups with boundary equal to the Pontryagin surfaces $\Pi_p$.  

Some other examples of boundaries come from a family of spaces defined as inverse limits of certain systems of connected sums of manifolds, introduced by 
Jakobsche~\cite{jakobsche}.  Such spaces are now known as \emph{trees of manifolds}.  
Fischer~\cite{fischer} showed that when the nerve of a right-angled Coxeter system is a connected, closed, orientable PL manifold, the boundary of the Davis complex is a tree of manifolds.  This was later corrected and extended to the non-orientable case by \'{S}wi\c{a}tkowski~\cite{swiat3}.  

Fischer's theorem 
did not automatically yield hyperbolic groups.  In~\cite{ps} Przytycki and \'{S}wi\c{a}tkowski construct flag-no-square triangulations of 3-dimensional PL manifolds.  This together with Fischer's result and Moussong's theorem (Theorem~\ref{moussong} above) 
yields hyperbolic right-angled Coxeter groups whose boundaries are trees of manifolds.  
Przytycki and \'{S}wi\c{a}tkowski also point out that Dranishikov's method together with  Fischer's result 
can be used to show that the Pontryagin sphere (which is different from the Pontryagin surfaces $\Pi_p$) occurs as the boundary of certain right-angled Coxeter groups.  (See Remarks 3.6 and 4.4(1) of~\cite{ps}.) 

Expanding to the non-hyperbolic case, in addition to the examples of Fischer mentioned above, \'{S}wi\c{a}tkowski~\cite{swiat2} gave sufficient conditions on the nerve of a right-angled Coxeter system which guarantee that the corresponding boundaries are $n$-dimensional Sierpinski compacta.  Roughly speaking, the criterion says that Sierpinski compacta occur if the nerve is an \emph{$(n+1)$-sphere with holes} 
(see~\cite[Definition 1.2]{swiat2}) such that the boundaries of the holes are ``nice'' and the holes are ``sufficiently far apart''.   

We end this section by mentioning work of \'{S}wi\c{a}tkowski~\cite{swiat4} which proposes a framework for exploring boundaries of groups which is more general than trees of manifolds.  
In~\cite{swiat4}, 
he introduces and develops the notion of trees of metric compacta and indicates that such spaces are potential candidates to be boundaries of groups.  

\subsection{Additional notes}
The results above focus on topological aspects of visual boundaries of (the Davis complexes of) right-angled Coxeter groups.  The interested reader may wish to consult the work of Hosaka (See \cite{hosaka} and the references therein) for some finer aspects of the actions of right-angled Coxeter groups on such boundaries.  

By considering the subset of points of the visual boundary of a CAT(0) group which come from the most 
``hyperbolic like'' geodesics, Charney and Sultan~\cite{charney-sultan} define \emph{contracting boundaries} which are invariant up to quasi-isometry.  These are discussed in Section~\ref{contracting}. 
For relatively hyperbolic groups, Bowditch~\cite{bowditch-relhyp} introduced another boundary, \emph{
the Bowditch boundary}, which generalizes the Gromov boundary of a hyperbolic group and is also a quasi-isometry invariant under some assumptions.  This is discussed in 
Section~\ref{bowditch}.

\section{Quasi-isometry classification}\label{sec:qi}

A central question in geometric group theory is that of classifying finitely generated groups up to quasi-isometry, and of understanding to what extent classes of groups are quasi-isometrically
rigid.  
A class is \emph{quasi-isometrically rigid} if any group quasi-isometric to a group in the class is virtually in the class. Such questions have been extensively explored for a number of classes of groups, for example lattices in semi-simple Lie groups, mapping class groups, 3-manifold groups and right-angled Artin groups. 

For right-angled Coxeter groups, one would like to obtain a quasi-isometry classification in terms of the defining graph $\Gamma$.   
We have already seen that right-angled Coxeter groups exhibit a variety of different (non-quasi-isometric) geometries (eg. hyperbolic, relatively hyperbolic, thick), and understanding the quasi-isometry classification within these classes is a challenging problem which will require an assortment of techniques, and very likely some new quasi-isometry invariants.   Below we outline the progress that has been made on this question so far.  
\subsection{Ends}

The number of ends of a space is a quasi-isometry invariant.  
It is well known that any finitely generated group has 0, 1, 2 or infinitely many ends~\cite{hopf}.  

The 0- and 2-ended cases are easy to understand.  All 0-ended groups are quasi-isometric to the trivial group, and all 2-ended groups are quasi-isometric to $\Z$.   Lemma~\ref{finite} and 
Theorem~\ref{two-ended} characterize the 0- and 2-ended right-angled Coxeter groups in terms of the defining graph.

For the case of infinitely many ends,
Theorem~\ref{one-ended} yields a visual 
tree-of-groups decomposition of the right-angled Coxeter group.  
 Now Theorems 0.3 and~0.4 of Papasoglu--Whyte~\cite{papa-whyte} imply that two right-angled Coxeter groups with infinitely many ends are quasi-isometric if and only if their 
tree-of-groups decompositions have the same set of quasi-isometry types of 1-ended vertex groups.  

It follows that classifying the one-ended right-angled Coxeter groups up to quasi-isometry is the key to understanding the quasi-isometry classification, and for the remainder of this section, we will restrict to 
one-ended case. This can be detected from $\G$ using Corollary~\ref{cor:one-ended}.

\subsection{Divergence and thickness}\label{divthick}
Given a pair of geodesic rays emanating from a basepoint, their divergence measures, as a function of $r$, the length of a shortest path which connects their time-$r$ points.  The divergence of a group is a quasi-isometry invariant which, roughly speaking, captures the largest possible divergence function exhibited by any pair of geodesics in a Cayley graph for the group.  See~\cite{behr-hag-sis}, \cite{dt1} or \cite{levco1} for a precise definition. 

A closely related quasi-isometry invariant is thickness, which was introduced in~\cite{behr-dru-mos}.  Thick spaces are defined inductively, with the base level, i.e.,~the spaces which are thick of order 0, corresponding to spaces with linear divergence.   Then intuitively, a space is \emph{thick of order $n$}, if it is a ``network'' of subsets, each of which is thick or order
 $n-1$,  with the property that any pair of points in the space can be connected by a chain of subsets which are thick of order $n-1$, such that successive subsets in the chain have infinite diameter intersection.  A refinement of this notion called \emph{strong algebraic thickness} which is better suited to finitely generated groups was introduced
in~\cite{behrstock-drutu}.   The structure of a network in the definition of thickness is well suited to obtaining upper bounds on divergence, and indeed it is proved in~\cite{behrstock-drutu} that any group which is strongly algebraically thick of order $n$ has divergence at most $x^{n+1}$.

By the Behrstock--Hagen--Sisto theorem (Theorem~\ref{relhypthick} above), 
a right-angled Coxeter group is either relatively hyperbolic or thick. 
Sisto showed that if a group is relatively hyperbolic, then its divergence is 
exponential~\cite[Theorem~1.3]{sisto}.  On the other hand it is well known that the divergence of a one-ended finitely presented group is at most exponential. (See~\cite[Lemma~6.15]{sisto}.)  Thus:
\begin{thm}\cite{behr-hag-sis}
The divergence of a right-angled Coxeter group is either exponential (if the group is relatively hyperbolic) or bounded above by a polynomial (if the group is thick).  
\end{thm}

In the thick case, the author and Thomas showed:
\begin{thm}\cite{dt1}\label{divxd}
For every positive integer $d$, there is a right-angled Coxeter group with divergence $x^d$.
\end{thm}

From the point of view of classifying right-angled Coxeter groups up to quasi-isometry, one would then like to address the following question.

\begin{qn}
Characterize the graphs $\G$ for which the divergence of $W_\G$ is $x^d$.  
\end{qn}

For linear divergence, one has the following. (The equivalence of (i) and (iii) below was proved 
in~\cite{dt1} for $\G$ triangle-free, while the general proof appears   in~\cite{behr-hag-sis})

\medskip
\begin{thm}~\cite{dt1, behr-hag-sis} 
For right-angled Coxeter groups $W_\G$, the following are equivalent: 
(i) the divergence of $W_\G$ is linear; (ii) $W_\G$ is thick of order 0; and (iii) $W_\G$ splits as a product or equivalently, $\G$ splits as a non-trivial join. 
\end{thm}
 
In order to characterize the quadratic divergence case, the author and Thomas defined a graph theoretic property called $\mathcal{CFS}$ on triangle-free graphs in~\cite{dt1}, which was generalized to all graphs in~\cite{bfhs}.
Roughly speaking, the $\mathcal{CFS}$ property says that the graph can be built out of squares in such a way that every two squares are connected by a ``chain'' of squares, where successive squares in a chain intersect along a diagonal.  (See~\cite{bfhs} for the definition in the general case.)  This condition guarantees that the 
corresponding Davis complex is built out of ``chains'' of flats, where successive flats in a chain have infinite diameter intersections.  As a result, the group is thick of order 1, and has quadratic divergence.   
In fact, the following is true:

\begin{thm}\cite{dt1,levco1} \label{racg-quad}
A right-angled Coxeter group $W_\G$ has quadratic divergence if and only if $\G$ is $\mathcal{CFS}$.  If $\G$ is not $\mathcal{CFS}$
then the divergence of $W_\G$ is at least cubic.  
\end{thm}
This was proved in~\cite{dt1} in the 2-dimensional case (when $\G$ is triangle free) and was generalized to all dimensions by Levcovitz in~\cite{levco1}.  The challenge here is to show that when the $\mathcal {CFS}$ property fails, the group actually has cubic divergence.  
This is done in~\cite{dt1} by analyzing certain van Kampen diagrams.  
In~\cite{levco1}, Levcovitz introduces a useful tool for obtaining such lower bounds on divergence in CAT(0) cube complex groups, called the \emph{hyperplane divergence function}.   


A full characterization of the right-angled Coxeter groups which have divergence $x^d$ for $d>2$ is still open, but Levcovitz makes significant progress in~\cite{levco1, levco2}. 
In~\cite{levco1} he defines a \emph{rank $n$ pair} to be a pair of non-adjacent vertices satisfying a certain inductive criterion.
Roughly, a non-adjacent pair is rank 1 if it is not contained in any square, and rank $n$ if 
the link of at least one of the vertices in the pair 
has the property that every pair of vertices is a rank $n-1$ pair.  
He shows that  the existence of such a pair implies a lower bound of $x^{n+1}$ on the divergence.  

In~\cite{levco2}, he defines the 
\emph{hypergraph index} (which for a thick group is always a non-negative integer), and proves that it is a quasi-isometry invariant  
for  2-dimensional right-angled Coxeter groups.  
Furthermore, he proves that having hypergraph index 
$n$ implies that the group is thick of order at most $n$ and has divergence at most $x^{n+1}$.  
Combining these two ideas, he obtains the following theorem, which in particular  recovers the examples of Theorem~\ref{divxd} above:

\medskip
\begin{thm}~\cite{levco1, levco2}
If $\Gamma$ contains a ``rank $n$ pair'' and has hypergraph index $n$, then $W_{\Gamma}$ has divergence $x^{n+1}$.
\end{thm}

The advantage of the hypergraph index over thickness and divergence is that it is defined by an algorithmic graph theoretic procedure, and is therefore always computable.  
\begin{conj} \cite{levco2}
The notions of hypergraph index $n$, thickness of order $n$ and divergence $x^{n+1}$ are equivalent.
\end{conj}
Levcovitz points out that the conjecture is true in the cases $n=0$ and $n=1$.

%
%

 \subsubsection{Random right-angled Coxeter groups, thickness, and relative hyperbolicity}\label{random-rh}
In~\cite[Section 3]{behr-hag-sis}, 
Behrstock, Hagen, and Sisto consider random right-angled Coxeter groups in the context of thickness and relative hyperbolicity, continuing the study begun by Charney--Farber.  
 (See Section~\ref{random} above).  In this setting, the probability function $p(n)$ can be thought of as an analog of the density in Gromov's density model.  
In~\cite[Figure 5]{behr-hag-sis}, the authors give a conjectural picture of threshold intervals for 
$p(n)$ which should correspond to relative hyperbolicity or thickness of various orders.  Roughly, 
the random right-angled Coxeter group should be relatively hyperbolic at low densities. Then there should be a sequence of non-overlapping threshold intervals of densities above that corresponding to groups that are thick of various orders, and finally, the groups should be virtually cyclic or finite at the highest densities.  
In~\cite[Section 3]{behr-hag-sis}, they prove a number of results which make progress towards this conjectural picture, including Theorem~3.4 which obtains a (low density) threshold for the random group to be relatively hyperbolic, and Theorem~3.9, which obtains (high density) thresholds for the 
random group to be finite, virtually cyclic and thick of order 0.  

This theme is continued in~\cite{bfhs}, 
where Behrstock, Falgas-Ravry, Hagen and Susse obtain a threshold interval for a random group to be thick of order 1:
\begin{thm}\cite[Corollary 3.2]{bfhs}\label{random-cfs}
There exists a constant $C>0$ such that if $ \left(\frac{C\log n}{n}\right)^\frac 12 \le p(n) \le 1- \frac{(1+\epsilon)\log n}{n} $ for some $\epsilon > 0$,  
then the random right-angled Coxeter group defined by a graph in $G(n, p(n))$ is a.a.s. thick of order exactly 1 and has quadratic divergence. 
\end{thm}
The authors observe that if $\G$ is $\mathcal{CFS}$ and $\G$ does not split as a non-trivial join, then $W_\G$ is thick of order exactly 1.  They obtain threshold probability functions which guarantee that a random graph has these properties, and use these to prove the above theorem. 
Moreover, they also show: 

\medskip
\begin{thm}~\cite[Theorem 5.7]{bfhs}
If $p(n) \le \frac{1}{\sqrt n \log n}$, then $\G\in G(n, p(n))$ is a.a.s. not in $\mathcal{CFS}$. 
\end{thm}
Combining this with Theorem~\ref{racg-quad} above, one sees that if $p(n)$ is as in the theorem, then a right-angled Coxeter group defined by a graph in $G(n, p(n))$ a.a.s.~has at least cubic divergence.  Behrstock, Falgas-Ravry, Hagen and Susse conjecture that there is a sharp threshold function below which $\G$ is a.a.s.~not $\mathcal{CFS}$, and above which $\G$ is 
a.a.s.~$\mathcal{CFS}$.   
See~\cite[bottom of page 4 and Remark~5.8]{bfhs} for a discussion of this conjecture.

 \subsection{Hyperbolic groups}

It is well-known  that hyperbolicity is a quasi-isometry invariant.  
Some of the simplest examples of hyperbolic right-angled Coxeter groups are ones whose defining graphs are $n$-cycles $\G_n$, where $n\ge 5$.  The generators of $W_{\G_n}$ act on $\Hyp^2$ by 
reflections in a right-angled hyperbolic $n$-gon.  
It is well known to the experts that the only Coxeter groups which act 
properly and cocompactly on $\Hyp^2$ (i.e., the only \emph{cocompact Fuchsian Coxeter groups}) are 
direct products of hyperbolic polygon reflection groups with finite (possibly trivial) Coxeter groups.  An explicit proof is recorded in~\cite[Appendix A]{dt2}.   By the Milnor--\v{S}varc lemma, these are all quasi-isometric to each other.  
Moreover work of Tukia~\cite{tukia}, Gabai~\cite{gabai} and 
Casson--Jungreis~\cite{casson-jungreis} shows that cocompact Fuchsian groups form a rigid quasi-isometry class, i.e., any finitely generated group quasi-isometric to a cocompact Fuchsian group is cocompact Fuchsian.  A direct proof of this fact for 2-dimensional right-angled Coxeter groups is given in~\cite{dt2}. 
It follows that 
\begin{thm}\label{cycles}
For 2-dimensional right-angled Coxeter groups, being cocompact Fuchsian, being quasi-isometric to a cocompact Fuchsian group, and having an $n$-cycle with $n\ge 5$ as defining graph are all equivalent. 
\end{thm}

Papasoglu~\cite{papasoglu} showed that the property of splitting over a 2-ended subgroup is a quasi-isometry invariant for groups which are not cocompact Fuchsian.

 For one-ended hyperbolic groups which are not co-compact Fuchsian and which split over a 2-ended subgroup, Bowditch~\cite{bowditch} defined a JSJ tree which encodes all of the canonical splittings over 2-ended subgroups.  A key feature of this tree is that its vertices (of which there are three types) and edges are defined in terms of subsets of local cut points in the visual boundary. Then it is immediate that given any quasi-isometry between two groups, 
the corresponding boundary homeomorphism (guaranteed by Theorem~\ref{hyp-quasi}) preserves these features 
of the visual boundary, 
and there is a type-preserving isometry between the corresponding JSJ trees.  In other words, Bowditch's JSJ tree is a quasi-isometry invariant.

 We now focus on $2$-dimensional right-angled Coxeter groups, i.e., the ones with $\G$ triangle-free.  Applying Theorems~\ref{thm:stallings},~\ref{moussong},~\ref{cycles}, and~\ref{two-ended-splittings} to this case, one easily obtains assumptions on the graph $\G$ which guarantee that $W_\G$ is 
 2-dimensional and satisfies the hypotheses of Bowditch's theorem:
 \begin{assumptions}\label{assumptions}
 The graph $\G$ 
 \begin{enumerate}
 \item has no triangles ($W_\G$  is 2-dimensional);
 \item  is connected and has no separating vertices or edges ($W_\G$ is 1-ended);
 \item   has no squares ($W_\G$  is hyperbolic);
 \item  is not a cycle of length $\ge 5$ ($W_\G$ is not cocompact Fuchsian); and
 \item has a separating pair of non-adjacent vertices ($W_\G$ splits over a 2-ended subgroup).
 \end{enumerate}
 \end{assumptions}

 For such groups, the author and Thomas gave a description of the JSJ tree in terms of $\G$:
 
 \medskip
 \begin{thm}~\cite{dt2}\label{jsjtree}
 For $\G$ satisfying Assumptions~\ref{assumptions}, the Bowditch JSJ tree for $W_\G$ can be described visually in terms of $\G$.  The vertices and edges of $\G$ correspond to subsets of vertices of $\G$ satisfying explicit graph theoretic conditions.   
 \end{thm}

A detailed statement of this visual description is given in~\cite[Theorem~3.37]{dt2}.
From this description it is immediate that the Bowditch tree can be computed algorithmically for these groups. 

The three types of vertices in Bowditch's tree 
are called \emph{finite-order, quadratically hanging (or maximal hanging Fuchsian), and rigid (or stars)}.  
The author and Thomas show that the JSJ tree of $W_\G$ with $\G$ satisfying Assumptions~\ref{assumptions}
has rigid vertices if and only if $\G$ as a $K_4$ minor (i.e., a subgraph which is a subdivided copy of $K_4$).  This is useful because when there are no rigid vertices, the Davis complex can be thought of as being built out of ``fattened trees'' of uniformly bounded valence, glued along their boundary components.  Then an isomorphism between the JSJ trees of two groups together with this ``tree of fattened trees'' structure can be used to construct a quasi-isometry between the groups themselves.  Such a construction was first used by Berhstock--Neumann 
in~\cite{behrstock-neumann}, and then by Malone in~\cite{malone} and Cashen--Martin in~\cite{cashen-martin}.
 
In light of the above discussion, we restrict to the following class of groups. 
\begin{definition}\label{nok4}
Let $\mathcal G$ be the class of graphs satisfying Assumptions~\ref{assumptions} which have no $K_4$ minors.  Let $\mathcal W_{\mathcal G}$ denote the class of Coxeter groups defined by graphs in $\mathcal G$.  
\end{definition}
For the class $\mathcal W_{\mathcal G}$  the author and Thomas show:

\medskip
\begin{thm}~\cite{dt2}\label{wg}
The Bowditch JSJ tree is a complete quasi-isometry invariant for $\mathcal W_{\mathcal G}$.
More precisely, two groups in $\mathcal W_{\mathcal G}$ are quasi-isometric if and only if there is a type preserving isomorphism between their JSJ trees.  

It follows that the quasi-isometry problem is decidable for $\mathcal W_{\mathcal G}$.  Moreover, the set of quasi-isometry classes of $\mathcal W_{\mathcal G}$  is recursively enumerable.
\end{thm}

This leaves a number of open questions:
\begin{qn}
What happens when there are rigid vertices? Can the ``no $K_4$ minor'' condition implicit in
Theorem~\ref{wg} 
be removed?
\end{qn}

In the trees obtained in Theorem~\ref{jsjtree}, the rigid vertices are all virtually free.  
Consequently, one can use techniques developed by Cashen and Macura~\cite{cashen-macura} to handle quasi-siometries preserving line patterns in free groups to analyze this situation.   Using this, 
Cashen, the author, and Thomas  show  (in an appendix to~\cite{dt2})
 that the situation is different when there are rigid vertices. 

\medskip
\begin{thm}~\cite[Appendix B]{dt2}
Let $W_n$ be the right-angled Coxeter group defined by a ``sufficiently subdivided'' copy of $K_n$
for $n\ge 4$, where $K_n$ is the complete graph on $n$ vertices.  Then the JSJ trees of $W_n$ with $n\ge4$  are all isomorphic.  
However, $W_m$ is quasi-isometric to $W_n$ if and only if $m=n$.
\end{thm}

Thus the Bowditch JSJ tree is not a complete invariant in the presence of rigid vertices, and in particular, the ``no $K_4$'' assumption cannot be removed.   Nevertheless, the techniques 
of~\cite{cashen-macura} and~\cite{cashen-martin} can probably be used to say more about the case in which the JSJ trees have rigid vertices.  

\begin{qn}
What happens when $\G$ is not triangle-free (so that the Davis complex of $W_\G$ is not 2-dimensional)?
\end{qn}
It should be possible to adapt the methods used in~\cite{dt2} to obtain a visual description of the Bowditch JSJ tree for the 2-dimensional case to higher dimensions, but the statements of the characterizations of the various types of vertices would necessarily be more complicated.  A discussion of this is included in the introduction 
to~\cite{dt2}.

\begin{qn}
What happens when the hyperbolic groups do not split over a 2-ended group?
\end{qn}
Any such groups have trivial JSJ tree, and the method of~\cite{dt2} does not apply. 
Theorem~\ref{two-ended-splittings} gives a characterization of the graphs which define 
the right-angled Coxeter groups which do not split over 2-ended subgroups.
There are many such graphs which define hyperbolic groups; an explicit example is the 
Petersen graph.    When the graph is triangle-free Theorem~\ref{thm:kk} tells us that the boundary of the group must be homeomorphic to the Sierpinski carpet or the Menger curve.  Thus one needs more refined quasi-isometry invariants here.

\subsection{Relatively hyperbolic groups}\label{relhyp}
The theorems of Caprace and Behrstock--Hagen--Sisto (Theorems~\ref{caprace} 
and~\ref{relhypthick} above) can be used to identify the 
right-angled Coxeter groups that are relatively hyperbolic, and also minimal peripheral structures for them.  At the moment we are very far from a full classification of the relatively hyperbolic right-angled Coxeter groups, 
but recently, there has been a burst in activity on exploiting known results on quasi-isometries between relatively hyperbolic groups (eg.~from~\cite{behr-dru-mos, drutu-sapir, groff}) and obtaining 
new invariants to distinguish particular (families of) relatively hyperbolic right-angled Coxeter groups.  We describe these results here.  
\subsubsection{Generalized theta graphs}
 In a precursor to~\cite{dt2}, the author and Thomas introduced generalized theta graphs to generalize the family of right-angled Coxeter groups studied in~\cite{crisp-paoluzzi}.  These are examples of graphs in the family ${\mathcal G}$ from 
Definition~\ref{nok4}, and are defined as follows.

\begin{definition}\label{defgentheta}
A \emph{generalized theta graph} is a graph with two distinguished vertices of valence $k\ge 3$, with $k$ branches 
connecting them. 
The $i$th branch has 
 $n_i$ valence 2 vertices and $n_i+1$ edges, with $n_k \ge \dots n_1 \ge 1$.
The \emph{linear degree} $l$ is the cardinality of $\{n_i \,|\, n_i=1\}$, and the \emph{hyperbolic degree} is $k-l$.
\end{definition}
It is immediate from Theorem~\ref{moussong} that a group defined by a generalized theta graph is 
hyperbolic if and only if the linear degree is at most 1.  
When the linear degree $l$ is at least 2, 
Theorem~\ref{caprace} shows that the group is hyperbolic relative to the special subgroup generated by the two vertices of valence 2 together with the~$l$ vertices along the paths with $n_i=1$.  
Hruska--Stark--Tran~\cite{hru-sta-tra} use the work of Dru\c{t}u--Sapir~\cite{drutu-sapir} to classify the relatively hyperbolic groups in this family up to quasi-isometry.  
Combining the results of~\cite{dt2} classifying the hyperbolic ones 
with the results of~\cite{hru-sta-tra} one has the following theorem:

\medskip
\begin{thm}\cite{dt2, hru-sta-tra} \label{gen-theta-qi}
Let $W$ and $W'$ be right-angled Coxeter groups defined by generalized theta graphs with linear degrees $l$ and $l'$ and 
hyperbolic degrees $h$ and $h'$ respectively.  Assume $l \le l'$.  Then $W$ and $W'$ are quasi-isometric if and only if one of the following holds. 
\begin{enumerate}
\item $l=l' \le 1$ and $h=h'$
\item $l=0, l'=1$ and $h= 2(h'-1)$
\item $l=l'=2$ and $h,h' \ge 1$
\item $l,l'\ge 3$ and $h, h' \ge 1$
\item $l,l'\ge 3$ and $h= h'=0$
\end{enumerate}
\end{thm}

One interesting consequence of this description is that there are infinitely many quasi-isometry classes among the hyperbolic groups in the above family (for example when $l=l'=0$, the hyperbolic degree determines the quasi-isometry class) but there are only three quasi-isometry classes among the relatively hyperbolic ones (i.e., the ones coming from the last three conditions in the theorem). 

\subsubsection{Divergence spectrum}
The divergence of a space essentially captures the larg\-est possible divergence function of any pair of geodesics in  the space.  Thus one loses all of the information about pairs of geodesics that diverge at a lower rate.  Charney proposed considering the whole spectrum of divergences seen in a group as a way of potentially capturing more information.  
Tran gave a precise definition of the divergence spectrum in~\cite{tran1}, proved it is a quasi-isometry invariant, and used it to distinguish a family of right-angled Coxeter groups. (See Figure 1 in~\cite{tran1} for a picture of the defining graphs of groups in this family.)  He points out that this particular family consists of relatively hyperbolic groups, and could also have been distinguished using existing results for relatively hyperbolic groups.  However, it is conceivable that such an approach could help distinguishing other groups.

\subsubsection{Bowditch Boundary}\label{bowditch}
Bowditch~\cite{bowditch-relhyp} introduced a boundary for relatively hyperbolic groups which generalizes the Gromov boundary of hyperbolic groups as well as the limit set of a geometrically finite Kleinian group.  There is an equivalent definition of relatively hyperbolic groups using a certain \emph{cusped space} which is hyperbolic, instead of the coned-off Cayley graph, and the Bowditch boundary is the Gromov boundary of this cusped space.  This is shown to be a quasi-isometry invariant in~\cite{groff}.  In~\cite{haul-ngu-tra} Haulmark--Nguyen--Tran provide more information on the structure preserved by a homeomorphism of Bowditch boundaries which is induced by a quasi-isometry of groups.  

In~\cite{dt2}, which concerns hyperbolic right-angled Coxeter groups, the visual construction of the JSJ tree is obtained by relating separation properties of  local cut points in the visual boundary to separation properties of  vertices in the defining graph.  The tree is then used as a quasi-isometry invariant.  In a similar manner, Haulmark--Nguyen--Tran relate cut points and \emph{non-parabolic cut pairs} in the Bowditch boundary of a relatively hyperbolic group to properties of the defining graph (see Theorems 1.1 and 1.3 of~\cite{haul-ngu-tra}).  This provides them with new quasi-isometry invariants, which they use to distinguish particular right-angled Coxeter groups
(see Examples 3.7 and 3.10 of~\cite{haul-ngu-tra}).  

Furthermore, they obtain a theorem similar to \'{S}wi\c{a}tkowski's result characterizing the right-angled Coxeter groups with Sierpinski carpet boundary (Theorem~\ref{racg-sierp} above): they 
give necessary and sufficient conditions for a  relatively hyperbolic right-angled Coxeter group with planar defining graph to have Bowditch boundary equal to $S^2$ or equal to the Sierpinski carpet.  
This allows them (in Example~4.15 of~\cite{haul-ngu-tra}) to distinguish (up to quasi-isometry) a pair of specific right-angled Coxeter groups which both have Sierpinski carpet 
visual boundary, but with distinct Bowditch boundaries (equal to $S^2$ in one case and the Sierpinski carpet in the other).  They also use the Bowditch boundary to show (in Lemma~5.7) that two right-angled Coxeter groups whose visual boundaries are both Menger curves are not quasi-isometric.

\subsection{Other leads for investigating quasi-isometry classification}
Here we mention a few other quasi-isometry invariants which have the potential to be useful for distinguishing right-angled Coxeter groups.

\subsubsection{Contracting boundaries}\label{contracting}
Recall that the boundary of a CAT(0) group is not a quasi-isometry invariant (and is not even well-defined).
Charney and Sultan~\cite{charney-sultan} defined \emph{contracting boundaries} for CAT(0) groups, which they proved to be quasi-isometry invariants.  (These were generalized to the setting of arbitrary proper geodesic metric spaces by Cordes~\cite{cordes}, and are now also called Morse boundaries.)  
Roughly speaking these boundaries pick up only the geodesic rays which demonstrate the most ``hyperbolic behavior''.  As a set, the contracting boundary is the collection of contracting (or equivalently Morse) geodesic rays.  When endowed with a certain direct limit topology, this becomes a quasi-isometry invariant. 
Charney and Sultan used the contracting boundary to distinguish two specific right-angled Coxeter groups.  
(See Figure 8 of~\cite{charney-sultan} for the defining graphs of these groups.)  

It follows from~\cite[Theorem F]{cordes-hume} (and is implicit in~\cite{charney-sultan})
that the contracting boundary of any right-angled Artin group is totally disconnected.   
Charney and Sisto asked if the same is true for thick right-angled Coxeter groups.  
A rather surprising example of Behrstock answers this in the negative, and also resolves some other questions in geometric group theory.  
\medskip
\begin{thm}~\cite{behrstock}\label{ex-behrstock}
There exists a right-angled Coxeter group $W_\G$ such that 
\begin{enumerate}
\item $W_\G$ has quadratic divergence, or equivalently, by Theorem~\ref{racg-quad}, $\G$ is $\mathcal{CFS}$
\item $W_\G$ contains a closed hyperbolic surface subgroup which is \emph{quasi-geodesically stable}
\item $W_\G$ contains a topologically embedded circle
in its contracting boundary.
\end{enumerate}
\end{thm}
In particular, (3) implies that the contracting boundary of $W_\G$ is 
not totally disconnected, and (2) answers a question of Taylor about quasi-geodesically stable subgroups.  
See Figure~1 of~\cite{behrstock} for a picture of a $\G$ which defines an example $W_\G$
as in Theorem~\ref{ex-behrstock} (there are endless variations, as will be evident from the following description).  It is implicit in the definition of the $\mathcal{CFS}$ property  that when one adds edges between existing vertices of a $\mathcal{CFS}$ graph, the resulting graph remains $\mathcal{CFS}$.  To construct his example, Behrstock adds edges to a 
$\mathcal{CFS}$ graph so that the new edges form an induced cycle of length $\ge 5$.  The 
quasi-geodesically stable surface subgroup is then a finite-index subgroup of the special subgroup of $W_\G$ corresponding to this cycle, and the topologically embedded circle in the 
contracting boundary is the boundary of this special subgroup.

Cordes and Hume~\cite{cordes-hume} have defined two quasi-isometry invariant notions of dimension for Morse boundaries, the stable dimension and the Morse capacity dimension, which may also prove useful.  

\subsubsection{Exploiting the hierarchically hyperbolic structure}
An important step in many quasi-isometric rigidity and classification results is showing that 
quasi-isometries take maximal dimensional flats close to flats.  In a recent paper, 
Behrstock--Hagen--Sisto~\cite{behr-hag-sis2} show that any maximal dimensional quasiflat in a hierarchically hyperbolic group lies within a finite distance of a union of standard orthants (under some mild conditions).  Moreover, they show that any quasi-isometry between hierarchically hyperbolic spaces induces a quasi-isometry between certain factored spaces, which are simpler than hierarchically hyperbolic spaces.  They outline a strategy for using these results to prove quasi-isometric rigidity and classification results for right-angled Coxeter and Artin groups.

\section{Commensurability questions}\label{sec:comm}

 Two groups are \emph{(abstractly) commensurable} if they have isomorphic finite-index subgroups.  
From this it is immediate that commensurability implies quasi-isometry.   In some classes of groups commensurability turns out to be equivalent to quasi-isometry  
(eg.~non-uniform lattices of $PSL_2(\mathbb C)$ are quasi-isometric if and only if they are commensurable~\cite{schwartz}), but this is not true in general (eg.~uniform lattices in  $PSL_2(\mathbb C)$ are all quasi-isometric but are not all commensurable~\cite{margulis}).  
\subsection{Commensurability classification in right-angled Coxeter groups}  
First consider the class of right-angled Coxeter groups defined by $n$-cycles with $n\ge 5$.  
Any such group is generated by reflections in the sides of a right-angled $n$-gon in $\Hyp^2$, which is finitely covered by a closed hyperbolic surface, and this in turn finitely covers the closed surface of genus two.  It follows that all groups in this class are commensurable, and therefore quasi-isometric to each other.  In particular, commensurability and quasi-isometry are equivalent for this class.

However, the situation is very different for right-angled Coxeter groups defined by generalized theta graphs. (See Definition~\ref{defgentheta}.)   For these groups the author, Stark and Thomas  
obtain certain Euler characteristic vectors as commensurability invariants~\cite{dst}.  Every right-angled Coxeter group $W_\G$ has a well-defined Euler characteristic, which is equal to the Euler characteristic of the orbicomplex 
obtained by taking the quotient of the Davis complex by the action of $W_\G$, but it can also be written as a formula involving the cardinalities of the vertex and edge sets of $\G$. (See~\cite[Definition 1.4]{dst}.)   Given a generalized theta graph with $k$ branches, there is an associated to it an \emph{Euler characteristic vector} in $\Q^n$ whose $i$th entry is the Euler characteristic of the special subgroup defined by the $i$th branch.   
Then one has:

\begin{thm}\cite{crisp-paoluzzi, dst}\label{gentheta}
Let $W$ and $W'$ be right-angled Coxeter groups defined by generalized theta graphs with linear degrees $l, l'\le1$ respectively.  The Euler characteristic vectors associated to $W$ and $W'$ determine whether or not $W$ and $W'$ are commensurable.  In particular, when $l=l'=0$, 
the groups $W$ and $W'$ are commensurable if and only if their Euler characteristics are commensurable.  
\end{thm}
Here two vectors are said to be commensurable if they have integer multiples that 
are equal.  The case of generalized theta graphs with three branches and linear degree 1 was done by Crisp--Paoluzzi in~\cite{crisp-paoluzzi}. Their result was generalized to all hyperbolic right-angled Coxeter groups defined by generalized theta graphs in~\cite{dst}.   Combining 
Theorems~\ref{gen-theta-qi} and~\ref{gentheta}, one has that there are infinitely many commensurablity classes in each quasi-isometry class containing a right-angled Coxeter group defined by a generalized theta graph.  For example, if $l=l'=0$ and $W$ and $W'$ have the same hyperbolic degree, then $W$ and $W'$ are quasi-isometric but not necessarily commensurable.

When the generalized theta graph defines a non-hyperbolic group, i.e., when the linear degree is $\ge 2$, then Hruska--Stark--Tran show: 
\begin{thm}\cite{hru-sta-tra}\label{thm:hst}
If the two groups have commensurable Euler characteristic vectors, then the groups are commensurable. 
\end{thm}
They remark that the following is still open, 
but indicate that it reduces to the case when the two groups in $W$ have the same linear degree.  
\begin{qn}
Is the converse of Theorem~\ref{thm:hst} true?
\end{qn}

Theorem~\ref{gentheta} fits into a larger commensurability classification program where progress is conceivable, namely the setting of the groups 
$\mathcal W_{\mathcal G}$ from Definition~\ref{nok4}: 

\begin{qn}
Classify $\mathcal W_{\mathcal G}$ up to commensurability.  
\end{qn}

This is a natural family to consider in light of Theorem~\ref{wg}, which enables one to completely understand the quasi-isometry classification of groups in this family.  One would then like to find complete commensurability invariants for the class of  groups which have the same JSJ tree.  

With this ultimate goal in mind, the author, Stark and Thomas~\cite{dst} considered the family of groups defined by \emph{cycles of generalized theta graphs}, which roughly speaking consist of 
(at least three) generalized theta graphs arranged in a circle by identifying them along their distinguished vertices.  (See~\cite[Definition 1.10]{dst} for a precise definition.)  
They associate to each group defined by such a graph a collection of Euler characteristic vectors, and show:
\begin{thm}
Let $W_\G$ and $W_{\G'}$ be two groups defined by cycles of generalized theta graphs.  There is a set of conditions which are both 
necessary and sufficient for $W_\G$ and $W'_{\G'}$ to be commensurable.  These conditions can be stated purely in terms of the associated Euler characteristic vectors and their commensurability classes, and therefore are determined $\G$ and $\G'$. 
\end{thm}
 
The proof of both directions make rather strong use of the fact that these generalized theta graphs are arranged in a cyclic pattern.  The unfortunate consequence of this is that generalizing this to all of $\mathcal G$ does not seem practical.  

\subsection{Commensurability between classes of groups}

Using commensurability to import information that is known for one class of groups to another class of groups is 
sometimes a fruitful strategy.  Most notably, tremendous milage has been obtained from Davis and Januszkiewicz' result~\cite{dav-jan} that every right-angled Artin group is commensurable to some right-angled Coxeter group.   For instance, this implies that right-angled Artin groups are linear, which is a key step of Agol's proof of the virtual Haken conjecture.

Here we highlight a few instances of commensurabilities between classes of groups where one of the classes  is right-angled Coxeter groups. 

\subsubsection{3-manifold groups}\label{3-manifolds}
One could ask which right-angled Coxeter groups are virtually 3-manifold groups. 
Davis--Okun~\cite{davis-okun} give a construction to show that if the nerve of $W_\G$ is planar, connected and triangle-free, then $W_\G$ acts properly on a contractible 3-manifold.
Droms~\cite{droms} showed that the cover of the Davis complex corresponding to the commutator subgroup of $W_\G$ (which is the complex $P_L$ from Definition~\ref{davis-cplx}) embeds in a 3-manifold if and only if $\G$ is planar. (Here $\G$ is not required to be triangl-free.) He used this to prove constructively that if $\G$ is planar, then $W_\G$ is virtually a 3-manifold group.  On the other hand he showed that the group defined by $K_{3,3}$ is not a 3-manifold group.  

\subsubsection{Geometric amalgams of free groups}
In~\cite{lafont}, Lafont introduced   
 \emph{surface amalgams}, which are spaces that are constructed by gluing surfaces with boundary along their boundary components in such a way that each boundary component is identified with at least two others and the resulting space has no free boundary components.  Thus the neighborhood of any point in the resulting space is homeomorphic to either an open disk in $\R^2$ or to $S \times I$, where $S$ is a star with $n\ge 3$ prongs and $I$ is an interval in $\R$. The fundamental groups of these spaces, called \emph{geometric amalgams of free groups}, are graphs of groups with edge groups all equal to $\Z$ and vertex groups that alternate between free groups and $\Z$'s.  
See~\cite{dst} or~\cite{lafont} for a precise definition.  

Surface amalgams play a crucial role in the proofs of the commensurability results of~\cite{crisp-paoluzzi, dst}.  Lafont~\cite{lafont} showed that geometric amalgams of free groups are topologically rigid, meaning that any isomorphism of their fundamental groups induces a 
homeomorphism of the surface amalgams.  
Crisp--Paoluzzi defined piecewise hyperbolic 
orbicomplexes for the groups they considered, and
then used topological rigidity to obtain their necessary conditions on commensurability by observing that the covers of their orbicomplexes which correspond to torsion-free finite-index subgroups of the right-angled Coxeter groups are necessarily surface amalgams.  

In~\cite{dst} the authors construct a piecewise  hyperbolic orbicomplex corresponding to each 
group in $\mathcal W_\mathcal G$ and again apply topological rigidity.  
(We remark that defining these new orbicomplexes is necessitated by the fact that topological rigidity fails for quotients of the Davis complex~\cite{stark}.)   These orbicomplexes admit particularly nice degree 16 covers which are surface amalgams.  In particular, every right-angled Coxeter group
$W_\G \in \mathcal W_\mathcal G$ is commensurable to a geometric amalgam of free groups, and the underlying graph of the latter is obtained by applying an easy transformation to $\G$.

It is also shown  in~\cite{dst} that every geometric amalgam of free groups whose defining graph is a tree is commensurable to a right-angled Coxeter group.  As a corollary, one has the commensurability classification of geometric amalgams of free groups whose defining graphs are trees of diameter at most 4.  

On the other hand, an example is given in~\cite{dst} of a geometric amalgam of free groups which is not quasi-isometric, hence not commensurable, to any right-angled Coxeter group. 
This leads to the question:
\begin{qn}
Which geometric amalgams of free groups are commensurable to right-angled Coxeter groups in $\mathcal 
W_{\mathcal G}$?
\end{qn}

Stark makes some progress on this in upcoming work~\cite{stark2}, where she studies the class $\mathcal C$ of 1-ended hyperbolic groups that are not cocompact Fuchsian and whose JSJ tree does not contain any rigid vertices.   A consequence of the work in 
in~\cite{dst} is that for groups in $\mathcal C$, the following properties are all equivalent:
being quasi-isometric to a right-angled Coxeter group, being quasi-isometric to a group which is generated by finite-order elements, and being quasi-isometric to a group with JSJ graph (i.e., the quotient of the JSJ tree by the group action) a tree.  In~\cite{stark2}, Stark produces a concrete graph theoretic criterion which is equivalent to the above properties.  Her criterion is in terms of certain graphs associated to degree refinements of JSJ trees; the latter are certain  finite matrices which encode the JSJ tree.  

Stark also shows that there exists a group in $\mathcal C$ which is quasi-isometric to a right-angled Coxeter group, but is not commensurable to any right-angled Coxeter group.  Thus the set of right-angled Coxeter groups within $\mathcal C$ is not quasi-isometrically rigid.

In the surface amalgams considered in~\cite{dst}, the gluing maps of the boundary components are 
all homeomorphisms.  Hruska-Stark-Tran~\cite{hru-sta-tra} explore what happens when the gluing maps are allowed to have degree $>1$.  Specifically the surface amalgams they consider 
$i$th boundary curve of the annulus is mapped to a simple closed curve in $S_i$ by a degree 
$m_i$ map.   
They show that each fundamental group in the class they consider is commensurable to a right-angled Coxeter group with defining graph a generalized theta graph.  They use this to carry out the quasi-isometry classification in their family and to conclude that the each group in the family of groups they consider is virtually a 3-manifold group, using the result of
Davis--Okun mentioned in Section~\ref{3-manifolds}.  (We remark that they also determine which groups in the family are actually 3-manifold groups.)  
\subsubsection{Right-angled Artin groups}
Davis--Januszkiewicz~\cite{dav-jan} showed that every right-angled Artin group is commensurable to
(and in fact is a finite index subgroup of) some right-angled Coxeter group.
It is natural to ask to what extent the converse is true: 

\begin{qn}\label{raag}
Obtain necessary and sufficient conditions on $\G$ such that 
$W_\G$ contains a right-angled Artin group of finite index. 
\end{qn}

It is quite easy to see that there are right-angled Coxeter groups which are not quasi-isometric, and hence not commensurable to any right-angled Artin group, using for example, the divergence function.  
As discussed in Section~\ref{divthick}, the functions $e^x$ or $x^d$ for $d \in \Z_+$ all occur as divergence functions of one-ended right-angled Coxeter groups, but the divergence of a one-ended right-angled Artin group is at most quadratic~\cite{behrstock-charney}.  Thus by Theorem~\ref{racg-quad} above, 
if $W_\G$ is one-ended and is commensurable to a right-angled Artin group, $\G$ must be
$\mathcal{CFS}$ or a non-trivial join.  Quite surprisingly, this is not enough.  Behrstock's example in Theorem~\ref{ex-behrstock} shows that there is a $\mathcal{CFS}$ right-angled Coxeter group which is not commensurable (or even quasi-isometric) to any right-angled Artin group (since its contracting boundary cannot be homeomorphic to the contracting boundary of a right-angled Artin group).

LaForge investigated Question~\ref{raag} in his PhD thesis~\cite{laforge}, resolving the special cases when $W_\G$ is hyperbolic or virtually free.  In the general setting, he focussed on 
\emph{visible subgroups} of right-angled Coxeter groups $W_\G$, whose 
generators are elements of $W_\G$ of the form $ab$, where $a$ and $b$ are not adjacent in $\G$. He identified certain graph theoretic criteria (the \emph{star cycle} and \emph{chain chord} conditions) which prevent visible subgroups from being right-angled Artin groups.  The answer to Question~\ref{raag} remains open in general.  

\bibliographystyle{amsalpha}
\bibliography{refs}
%

\end{document}